\documentclass[11pt,reqno]{amsart}
\usepackage{graphicx,color}
\usepackage{amssymb,amscd,latexsym, mathrsfs, epsfig}
\usepackage[all]{xy}

\newcommand\PP{\mathbb P}
\newcommand\C{\mathbb C}

\newcommand\Aut{\operatorname{Aut}}

\makeatletter \@addtoreset{equation}{section} \makeatother

\title[]{A note on the definition of K-stability}
\author{Jacopo Stoppa}
\date{}
\address{}
\email{}
\begin{document}
\begin{abstract} As recently pointed out by Li and Xu, the definition of K-stability, and the author's proof of K-stability for cscK manifolds without holomorphic vector fields, need to be altered slightly: the Donaldson-Futaki invariant is positive for all test configurations \emph{which are not trivial in codimension 2}.
\end{abstract}
\maketitle
Let $(X, L)$ be a polarized manifold and $(\mathcal{X}, \mathcal{L})$ denote a test configuration. The usual definition of K-stability (see e.g. \cite{rich} Theorem/Definition 3.9) requires that the Donaldson-Futaki invariant $F(\mathcal{X}, \mathcal{L})$ is \emph{strictly positive} for all \emph{nontrivial} test configurations. In particular if $F(\mathcal{X}, \mathcal{L}) = 0$ then $(\mathcal{X}, \mathcal{L})$ should be $\C^*$-equivariantly isomorphic to the trivial test configuration $X \times \C$ with the trivial $\C^*$-action. (It follows that if $(\mathcal{X}, \mathcal{L})$ is K-stable then $\Aut(X, L)$ is finite). As recently pointed out by Li and Xu (see \cite{lixu} Section 2.2), this is asking too much, and the definition must be slightly modified. This is because for \emph{every} $(X, L)$ one can in fact construct many nontrivial (but almost trivial) test configurations with vanishing Donaldson-Futaki invariant, e.g. by identifying two points on $X$ in a flat way (see \cite{lixu} Remark 4 (2)). The simplest explicit example is the usual flat projection of a twisted cubic in $\PP^3$ to a rational nodal cubic with an embedded point. As pointed out in \cite{lixu} Example 1 this can be seen as a test configuration with vanishing Donaldson-Futaki invariant, which is not trivial or induced by a holomorphic vector field on $\PP^1$. To get rid of holomorphic vector fields altogether we can repeat the same example with a high genus curve in $\PP^3$ which projects flatly to a plane curve with a number of nodes, carrying embedded points. All these test configurations have a common feature: the $\C^*$-action on the central fibre $\mathcal{X}_0$ is trivial in codimension 1 (\emph{even on the nilpotents}, so e.g. the degeneration of a conic to a double line is \emph{not} an example, and indeed it has positive Futaki).

\smallskip

\noindent \textbf{Definition 1.} A test configuration $(\mathcal{X}, \mathcal{L})$ is \emph{trivial in codimension 2} if $\mathcal{X}$ is $\C^*$-equivariantly isomorphic to the product test configuration $X \times \C$ (with trivial $\C^*$-action), away from a closed subscheme of codimension at least $2$.

\smallskip

Accordingly, we should amend the definition of K-stability.

\smallskip

\noindent \textbf{Definition 2.} A polarised variety $(X, L)$ is K-stable if $F(\mathcal{X}, \mathcal{L}) > 0$ for all test configurations \emph{which are not trivial in codimension 2}. 

\smallskip

Recall the main result of \cite{io} (refining the K-semistability proved by Donaldson):

\smallskip

\noindent \textbf{Theorem (\cite{io} Theorem 1.2).} \emph{A polarised variety $(X, L)$ for which the group $\Aut(X, L)$ is finite and which admits a constant scalar curvature K\"ahler metric in the class $c_1(L)$ is K-stable.}

\smallskip

The purpose of this note is to correct this statement and its proof slightly, according to the discussion above: in other words, for this statement to be true, we need to use the definition of K-stability above (avoiding test configurations which are trivial in codimension $2$).  

The point is that in \cite{io} Proposition 3.3 we need to assume that $(\mathcal{X}, \mathcal{L})$ in not trivial in codimension $2$. The oversight in the proof  concerns the ``degenerate case", p. 9. We follow the setup and notation explained there. Recall that we proved that in the degenerate case there is a finite map $\rho: \mathcal{X}_1 \to \mathcal{X}^{\rm red}_0$ from the general fibre to the reduced central fibre. The problem is that we neglected to consider the case when this map is \emph{generically injective}. This happens precisely when the test configuration $(\mathcal{X}, \mathcal{L})$ is trivial in codimension 2. Suppose then that this is not the case, i.e. that $\deg\rho > 1$. Then $\mathcal{X}_0$ is generically nonreduced. Therefore we can choose the section $x_{r + i}$ of $\mathcal{L}$ which appears in the argument so that its restriction to $\mathcal{X}^{\rm red}_0$ is \emph{generically nonzero}. (This would fail in the case when $\rho$ is generically injective and so $\mathcal{X}_0$ is generically reduced). Multiplying the sections of $\mathcal{L}^{k-1}_{0}|_{\mathcal{X}^{\rm red}_0}$ by $x_{r + i}$ gives the upper bound $(3.7)$ for the weight $w(k)$ (which again would fail if $x_{r + i}$ is not generically nonzero), and the rest of the argument goes through without changes, proving that $F(\mathcal{X}, \mathcal{L})$ is strictly positive.

\smallskip
 
\noindent\textbf{Remark.} It follows from Definition 2 that it is enough to test K-stability with respect to \emph{normal} test configurations. Suppose that $(\mathcal{X}, \mathcal{L})$ is a test configuration with nonpositive Futaki invariant, with $\mathcal{X}$ not normal, and not trivial in codimension $2$. Let us denote by $(\mathcal{X}', \mathcal{L}')$ the normalisation of $\mathcal{X}$ with the pullback line bundle. According to \cite{rich} Remark 5.2 the Donaldson-Futaki invariant decreases with normalisation, so $F(\mathcal{X}', \mathcal{L}') \leq 0$. Suppose by contradiction that $(\mathcal{X}', \mathcal{L}')$ is trivial in codimension $2$. Then the original test configuration $\mathcal{X}$ must be degenerate in the sense of \cite{io}, proof of Proposition 3.3: if not, the $\C^*$-action of $\mathcal{X}_0$ is nontrivial on a Zariski open subset, and so the induced action on $\mathcal{X}'_0$ cannot be trivial in codimension 1. But on the other hand we know $F(\mathcal{X}, \mathcal{L}) \leq 0$, and so we cannot have $\deg \rho > 1$ (otherwise the Donaldson-Futaki invariant would be strictly positive). So we have reached a contradiction (since $\deg \rho = 1$ implies that $\mathcal{X}$ is in fact trivial in codimension 2). So one could give an equivalent definition of K-stability (in the sense of Definition 2) by requiring that $F > 0$ \emph{for all normal test configurations except the trivial one}. This is the approach taken by Li an Xu (\cite{lixu}, Definition 3). As they observe, it seems very likely that all the results about K-stability in the literature work perfectly fine if one restricts to normal test configurations (or equivalently adopting Definition 2). An alternative point of view has been suggested by Sz\'ekelyhidi in \cite{gabor}: it seems very likely that it is enough to test K-stability with respect to all test configurations of positive $L^2$ norm.

\smallskip

\noindent\textbf{Acknowledgements.} The author is grateful to C. Li and C. Xu for sharing with him a preliminary version of their work \cite{lixu}.

\end{document}